\magnification=\magstep1
\input amstex
\documentstyle{amsppt}
\pagewidth{6.5truein}
\pageheight{9truein}
\loadbold

\def\R{{\bold R}}

\def\N{{\bold N}}

\def\phi{\roman{\varphi}}

\def\Reg{\mathop{\roman{Reg}}\nolimits}

\def\const{\mathop{\roman{const}}\nolimits}

\def\graph{\mathop{\roman{graph}}\nolimits}
\def\id{\mathop{\roman{id}}\nolimits}
\def\dim{\mathop{\roman{dim}}\nolimits}

\def\Im{\mathop{\roman{Im}}\nolimits}

\def\bdry{\mathop{\roman{bdry}}\nolimits}
\def\Reg{\mathop{\roman{Reg}}\nolimits}
\def\proj{\mathop{\roman{proj}}\nolimits}
\def\epsilon{\varepsilon}
\topmatter
\title 
Triangulation of the map of a $G$-manifold to its orbit space
\endtitle 
\author 
Mitsutaka MURAYAMA, Masahiro SHIOTA \endauthor 
\address
Department of Mathematics, Tokyo Institute of Technology, 
Meguro-ku, Tokyo 152-8551, Japan\newline
Graduate School of Mathematics, Nagoya University, Chikusa, Nagoya
464-8602, Japan
\endaddress
\abstract
Let $G$ be a Lie group and $M$ a smooth proper $G$-manifold. 
Let $\pi:M\to M/G$ denote the natural map to the orbit space. 
Then there exist a PL manifold $P$, a polyhedron $L$ and homeomorphisms 
$\tau:P\to M$ and $\sigma:M/G\to L$ such that $\sigma\circ\pi\circ\tau$ is PL. 
If $M$ and the $G$-action are of analytic class, we can choose subanalytic 
$\tau$ and then unique $P$ and $L$. 
\endabstract
\email
murayama\@math.titech.ac.jp, 
shiota\@math.nagoya-u.ac.jp
\endemail
\endtopmatter
\document
\head 1. Introduction \endhead
Let $G$ be a Lie group, which we regard as of class $C^\omega$. 
A $C^k$ {\it $G$-manifold} $M,\,k=1,...,\infty,\omega$, is a $C^k$ manifold 
with an action of $G$ on $M$ such that the map $G\times M\ni(g,x)\to g x\in M$ 
is of class $C^k$. 
A manifold means a manifold without boundary, though our arguments work also in 
the case of manifolds with boundary. 
A $C^k$ $G$-manifold $M$ is called {\it proper} if the map $G\times M\ni(g,x)
\to(x,g x)\in M^2$ is proper. 
Let $M/G$ and $\pi:M\to M/G$ denote the orbit space $\{Gx\,|\,x\in M\}$ with the 
quotient topology and the natural map respectively. 
A {\it triangulation} of a topological space $X$ is a pair of a polyhedron 
and its homeomorphism to $X$. 
A {\it triangulation} of a $C^0$ map between topological spaces $\phi:X\to Y$ is 
a pair of triangulations $\tau_X:P_X\to X$ of $X$ and $\tau_Y:P_Y\to Y$ of $Y$ 
such that $\tau^{-1}_Y\circ\phi\circ\tau_X:P_X\to P_Y$ is PL. 
A $C^k$ {\it triangulation} of a $C^k$ manifold $N$ is a pair of a PL manifold 
$P$ with its simplicial decomposition $K$ and a homeomorphism $\tau:P\to N$ 
such that $\tau|_\sigma$ for each $\sigma\in K$ is a $C^k$ diffeomorphism onto 
$\tau(\sigma)$. 
Existence of a $C^k$ triangulation and uniqueness of the PL manifold up to PL 
homeomorphism are well-known (Cairns-Whitehead). 
Moreover, triangulability of an orbit space is considered by many people (e.g., 
Matumoto-Shiota [M-S$_{1,2}$], Verona [V] and Yang [Y]). 
In this paper we show a triangulation of $\pi:M\to M/G$. \par 
  For the construction of a triangulation of $\pi$ we give the subanalytic 
structures to $\pi:M\to M/G$, proceed in the category of subanalytic sets and maps 
and apply its theory. 
(Gabrielov introduced the category.)  
A {\it subanalytic} set is a subset of a Euclidean space $\R^n$ of the form 
$\cup_i(\Im f_{i1}-\Im f_{i2})$, where $f_{i j}$ are a finite number of proper 
real analytic maps of real analytic manifolds into $\R^n$, and a {\it 
subanalytic} map is a continuous map between subanalytic sets with subanalytic 
graph. 
Examples of a subanalytic set and a subanalytic map are a polyhedron closedly 
included in a Euclidean space and a PL map between such polyhedra respectively. 
Note a $C^\omega$ submanifold of $\R^n$ and a $C^\omega$ map between 
subanalytic $C^\omega$ submanifolds of $\R^n$ are not necessarily subanalytic. 
A sufficient condition for them to be subanalytic is that the submanifolds are 
closedly included in $\R^n$. 
We always assume $G$ satisfies this condition though some of forthcoming 
analytic submanifolds do not. 
Triangulations of a subanalytic set and a subanalytic map are called {\it 
subanalytic} if the polyhedra and the homeomorphisms are subanalytic. 
We define also a {\it subanalytic} $C^k$ $G$-manifold $M$ by requiring $M$ and 
the map $G\times M\to M$ are subanalytic. 

\proclaim{Theorem (Triangulation of $\pi:M\to M/G$)}
Let $G$ be a Lie group, and $M$ a proper $C^k$ $G$-manifold, $k=1,...,\infty,\omega$. 
Then there exist a triangulation of $\pi:M\to M/G$ $(\tau:P\to M,\,\sigma:L\to M/
G)$ such that $P$ is the PL manifold of a $C^k$ triangulation of $M$. 
Moreover, if $M$ is a subanalytic $C^k$ $G$-manifold, we can choose subanalytic 
$\tau$ and then $L$ is unique up to PL homeomorphism. 
\endproclaim

  Here, even if $M$ is a proper $C^\omega$ $G$-manifold, we do not know whether 
$\sigma^{-1}\circ\pi\circ\tau:P\to L$ is unique, i.e., whether for another 
triangulation $(\tau':P'\to M,\,\sigma':L'\to M/G)$ there exist PL homeomorphisms 
$\phi:P\to P'$ and $\psi:L\to L'$ such that $\psi\circ\sigma^{-1}\circ\pi\circ\tau
=\sigma^{\prime-1}\circ\pi\circ\tau'\circ\phi$. 
Another open problem is whether we can choose a triangulation $(\tau:P\to M,\,\sigma
:L\to M/G)$ of $\pi$ so that for each element $g$ of $G$, the action $P\ni x\to\tau
^{-1}(g\tau(x))\in P$ is PL. \par
  The theorem holds true also in the case where $M$ is a proper $C^k$ $G$-manifold 
with boundary as follows. 
We give naturally a proper $C^k$ $G$-manifold structure to the double $DM$ of $M$, 
and consider the pair of $DM$ and $\partial M$. 
Then it suffices to generalize theorem to the following form. 
In theorem, let $M'$ be a proper $C^k$ $G$-submanifold of $M$, and assume $M'$ is 
closed in $M$. 
Then we can choose the triangulation $(\tau:P\to M,\,\sigma:L\to M/G)$ so that 
$\tau^{-1}(M')$ is a subpolyhedron of $P$, which is clear by the following proof. 
\head
2. Subanalytic sets
\endhead
  The main idea of proof is to apply the affirmative answer to Thom's conjecture 
that a proper Thom map is triangulable [S$_2$]. 
It is easy to see $\pi:M\to M/G$ can be a Thom map. 
Hence in the case of compact $G$, $\pi$ is proper and theorem follows. 
However, if $G$ is not compact, $\pi$ may be non-proper and a non-proper Thom map 
is not necessarily triangulable. 
We modify Thom's conjecture as shown below. 
For that we need to proceed in the subanalytic category. 
We prepare some terminology and facts (see [G-al] for Whitney stratifications and 
tube systems, and [H] and [S$_1$] for the subanalytic category). 
Elementary properties of subanalytic sets are that for subanalytic subsets $X$ and 
$Y$ of $\R^n$, $X\cap Y,\ X\cup Y,\ X\times Y,\ X-Y$ and $\overline X$ are 
subanalytic, $\overline X-X$ is of dimension smaller than $X$ when $X\not=\emptyset$, 
the family of connected components of $X$ is locally finite at each point of $\R^n$, 
each connected component is subanalytic and a subanalytic $C^\infty$ manifold is of 
class $C^\omega$. 
For a subanalytic set $X\subset\R^n$, let $\Reg X$ denote the subset consisting 
of points where the germ of $X$ is $C^\omega$ smooth and of maximal dimension. 
A{\it subanalytic} $C^\omega$ {\it stratification} $\{X_i\}$ of $X$ is a partition 
of $X$ into a finite number of subanalytic $C^\omega$ manifolds $X_i$. 
It is known that $\Reg X$ is subanalytic and $\dim(X-\Reg X)<\dim X$. 
Hence $\{\Reg X, \Reg(X-\Reg X),...\}$ is the ``canonical" subanalytic $C^\omega$ 
stratification of $X$. 
We say 2 subanalytic $C^\omega$ manifolds $X_1$ and $X_2$ in $\R^n$ satisfy the 
{\it Whitney condition} at a point $b$ of $X_2\cap\overline{X_1}$ if the following 
statement holds. 
Let $\{a_k\}$ and $\{b_k\}$ be sequences in $X_1$ and $X_2$, respectively, both 
converging to $b$ such that the sequence of the tangent spaces $\{T_{a_k}X_1\}$ 
converges to a subspace $T\subset\R^n$ in $G_{n,m}$---the Grassmannian of 
$m$-dimensional subspaces of $\R^n$, where $m=\dim X_1$, and the sequence of the 
lines $\{\overrightarrow{a_kb_k}\}$ converges to a line $L\subset\R^n$ in $G_{n,1
}$. 
Then $L\subset T$. 
A {\it Whitney} subanalytic $C^\omega$ stratification $\{X_i\}$ is the case where 
each pair of $X_i$ and $X_{i'}$ satisfy the Whitney condition at each point of 
$X_{i'}\cap\overline{X_i}$. 
We know also that for 2 disjoint subanalytic $C^\omega$ manifolds $X_1$ and $X_2$ 
with $X_2\subset\overline{X_1}-X_1$ the subset of $X_2$ where $X_1$ and $X_2$ 
satisfy the Whitney condition is subanalytic, that its complement in $X_2$ is of 
smaller dimension than $\dim X_2$ and that the closure of a subanalytic set in the 
ambient Euclidean space is subanalytic and of the same dimension. 
Hence we can construct the ``canonical" Whitney subanalytic $C^\omega$ 
stratification $\{X_i\}$ of $X$ as follows. 
Let $X_1$ be $\Reg X,\,X_2$ be the union of $A_1=\Reg(X-\Reg X)-\overline{\Reg X}$ 
and the subset $A_2$ of $\Reg(\Reg(X-\Reg X)\cap\overline{\Reg X})$ where $\Reg X$ 
and $\Reg(X-\Reg X)$ satisfy the Whitney condition if both $A_1$ and $A_2$ are of 
the same dimension and the set of larger dimension otherwise, and $X_3,...,$ be so 
on. \par
  For a subanalytic map $\phi:X\to Y$ from a $C^\omega$ manifold in $\R^n$ to a 
subanalytic set in $\R^n$, let $\Reg\phi$ denote the points of $X$ where the germ 
of $\phi$ is $C^\omega$ smooth and has locally the maximal Jacobian rank. 
Then $\Reg\phi$ is subanalytic and $\dim(X-\Reg\phi)
<\dim X$. 
In general, let $\phi:X\to Y$ be a subanalytic map between subanalytic sets in 
$\R^n$. 
A {\it subanalytic $C^\omega$ stratification} of $\phi$ is a pair of subanalytic 
$C^\omega$ stratifications $\{X_i\}$ of $X$ and $\{Y_j\}$ of $Y$ such that for 
each $i$, $\phi|_{X_i}$ is a $C^\omega$ submersion to some $Y_j$. 
We write as $\phi:\{X_i\}\to\{Y_j\}$ and call it a {\it subanalytic $C^\omega$ 
stratified map}. 
Note if $X$ and $Y$ are bounded in $\R^n$, $\Im\phi$ and $\phi^{-1}(A)$ for any 
subanalytic subset $A$ of $Y$ are subanalytic. 
Hence, if \newline
$(*)$ $\phi(B\cap X)$ and $\phi^{-1}(B)$ for each bounded set $B$ in $\R^n$ 
are bounded in $\R^n$, \newline
then by the above fact in the case of smooth $X$ there always exists the 
``canonical" subanalytic $C^\omega$ stratification $\phi:\{X_i\}\to\{Y_j\}$ of $\phi$. 
However, this is not the case without $(*)$ in general. 
For example, let $X=\N$, $Y=\R$ and $\phi$ defined so that $\phi(k)=1/k$ for 
$k\not=0$. 
Then $\phi$ does not admit a subanalytic $C^\omega$ stratification. 
We call a subanalytic $C^\omega$ stratified map $\phi:\{X_i\}\to\{Y_j\}$ a 
{\it Whitney} stratification if $\{X_i\},\,\{Y_j\}$ and $\{\graph\phi|_{X_i}\}$ 
are Whitney stratifications. 
Under the above condition $(*)$ it follows also that there exists the canonical 
Whitney subanalytic $C^\omega$ stratification of $\phi$.
A $C^\omega$ {\it function} on a subset $X\subset\R^n$ is the restriction to $X$ 
of a $C^\omega$ function defined on an open neighborhood of $X$ in $\R^n$. 
A $C^\omega$ {\it map} from a subset $X\subset\R^n$ to another $Y\subset\R^n$ is 
defined in the same way. 
Note if $X$ is a $C^\omega$ submanifold of $\R^n$, this definition is equivalent to 
that the map is analytic in the usual sense. 
If the underlying map $\phi:X\to Y$ of a $C^\omega$ stratified map $\phi:\{X_i\}
\to\{Y_j\}$ is of class $C^\omega$, the Whitney condition on $\{\graph\phi|_{X_i}
\}$ in the definition of a Whitney stratification is not necessary; and in our 
arguments we can replace $\phi:X\to Y$ with $\proj:\graph\phi\to Y$. 
Hence for simplicity we always consider maps with stratification are of class $C^
\omega$. \par 
  Let $\{A_\alpha\}$ and $\{B_\beta\}$ be families of subanalytic subsets of 
subanalytic sets $X,\,Y\subset\R^n$, respectively, locally finite at each point of 
$\R^n$. 
Then a subanalytic $C^\omega$ stratification $\{X_i\}$ of $X$ is {\it compatible 
with} $\{A_\alpha\}$ if each $A_\alpha$ is a union of some connected components of 
$X_i$'s. 
The canonical Whitney subanalytic $C^\omega$ stratification $\{X_i\}$ of $X$ 
compatible with $\{A_\alpha\}$ exists. 
Indeed, we define $X_1$ to be $\Reg X-\cup\{\overline{A_\alpha}\,|\,\dim A_\alpha<
\dim X\}-\cup\{\overline{A_\alpha}-\Reg A_\alpha\,|\,\dim A_\alpha=\dim X\}$, $X_2$ 
by considering $X-X_1,\,\{A_\alpha-X_1\}$ and the Whitney condition, and so on. 
We define naturally also a subanalytic $C^\omega$ stratification of a subanalytic 
map $\phi:X\to Y$ {\it compatible with} $\{A_\alpha\}$ and $\{B_\beta\}$. 
Under the above condition $(*)$ we can construct in the same way as above the 
canonical Whitney subanalytic $C^\omega$ stratification $\phi:\{X_j\}\to\{Y_k\}$ 
of $\phi$ compatible with $\{A_\alpha\}$ and $\{B_\beta\}$. \par
  A Whitney subanalytic $C^\omega$ stratified map $\phi:\{X_i\}\to\{Y_j\}$ is 
called a {\it Thom map} if the following condition is satisfied. 
Let $X_i$ and $X_{i'}$ be strata such that $\overline{X_i}\cap X_{i'}\not=\emptyset
$. 
If $\{a_k\}$ is a sequence of points of $X_i$ convergent to a point $b$ of $X_{i'}
$, and if the sequence of the tangent spaces $\{T_{a_k}(\phi|_{X_i})^{-1}(\phi(a_k
))\}$ converges to a space $T\subset\R^n$ in $G_{n,l},\,l=\dim(\phi|_{X_i})^{-1}(
\phi(a_k))$, then $T_b(\phi|_{X_{i'}})^{-1}(\phi(b))\subset T$. \par
  Let $m\ge 2$ be an integer. 
A {\it subanalytic $C^m$} (not $C^\omega$) {\it tube system} $\{T_j=(|T_j|,\pi_j,
\rho_j)\}$ for a Whitney subanalytic $C^\omega$ stratification $\{Y_j\}$ of a 
subanalytic set $Y\subset\R^n$ consists of one tube $T_j$ at each $Y_j$, where $\pi
_j:|T_j|\to Y_j$ is a subanalytic $C^m$ submersion of an open tubular neighborhood 
of $Y_j$ in $\R^n$ and $\rho_j$ is a non-negative subanalytic $C^m$ function on $|T_
j|$ such that $\rho_j^{-1}(0)=Y_j$ and each point $y$ of $Y_j$ is a unique and 
non-degenerate critical point of $\rho_j|_{\pi^{-1}_j(y)}$. 
We call a tube system $\{T_j\}$ {\it strongly controlled} if for each pair $j$ 
and $j'$ with $\dim Y_j<\dim Y_{j'}$, the following properties hold true\,:
$$
\gather
\pi_j\circ\pi_{j'}=\pi_j\quad\text{and}\quad\rho_j\circ\pi_{j'}=\rho_j\quad
\text{on}\ |T_j|\cap|T_{j'}|,\tag ct
\endgather
$$
and\newline
(sc) the map $(\pi_j,\rho_j)_{Y_{j'}\cap|T_j|}$ is a $C^m$ submersion into $Y_j
\times\R$. \newline
Note that (sc) follows from (ct) since $(\pi_j,\rho_j)|_{Y_{j'}\cap|T_j|}\circ\pi
_{j'}=(\pi_j,\rho_j)$ on $|T_j|\cap|T_{j'}|$, hence our definition of strongly 
controlledness coincides with that of controlledness in [G-al] and [S$_1$], and 
that any Whitney subanalytic $C^\omega$ stratification admits a strongly controlled 
subanalytic $C^m$ tube system [S$_1$, Lemma I.1.3]. 
(The reason why we consider (sc) will become clear. 
There does not necessarily exist a strongly controlled subanalytic $C^\infty(=C^
\omega)$ tube systems. 
Hence we need to consider class $C^m$.) 
Let $\phi:\{X_{i,j}\}_{j=1,...,k\atop i=1,...,l_j}\to\{Y_j\}_{j=1,...,k}$ be a 
Whitney subanalytic $C^\omega$ stratification of a subanalytic $C^\omega$ map $\phi
:X\to Y$ such that $\phi^{-1}(Y_j)=\cup_{i=1}^{l_j}X_{i,j}$ for each $j$, where $X,
\,Y\subset\R^n$. 
Let $\{T_j=(|T_j|,\pi_j,\rho_j)\}$ be a strongly controlled subanalytic $C^m$ 
tube system for $\{Y_j\}$, and let $\{T_{i,j}=(|T_{i,j}|,\pi_{i,j},\rho_{i,j})\}$ 
be a subanalytic $C^m$ tube system for $\{X_{i,j}\}$. 
We call $\{T_{i,j}\}$ {\it strongly controlled over} $\{T_j\}$ if the following 
conditions are satisfied. \newline
(sc1) For each $(i,j)$, $\tilde\phi(|T_{i,j}|)\subset|T_j|$ and $\phi\circ\pi_{i,j}
=\pi_j\circ\tilde\phi$ on $|T_{i,j}|$, where $\tilde\phi$ is a subanalytic $C^\omega$ 
extension of $\phi$ to a subanalytic open neighborhood of $X$ in $\R^n$. \newline
(sc2) For each $j$, $\{T_{i,j}\}_{i=1,...,l_j}$ is a strongly controlled tube 
system for $\{X_{i,j}\}_{i=1,...,l_j}$. 
(sc3) For any pair $(i,j)$ and $(i',j')$ with $\dim Y_j<\dim Y_{j'}$ and $\dim X_{
i,j}<\dim X_{i',j'}$, it holds that $\pi_{i,j}\circ\pi_{i',j'}=\pi_{i,j}$ on $|T_{
i,j}|\cap|T_{i',j'}|$ and $(\pi_{i,j},\phi)_{X_{i',j'}\cap|T_{i,j}|}$ is a $C^m$ 
submersion into the $C^m$ manifold $\{(x,y)\in X_{i,j}\times(Y_{j'}\cap|T_j|)\,|\,
\phi(x)=\pi_j(y)\}$. \par
  If the latter condition in (sc3) fails, $\{T_{i,j}\}$ is called {\it controlled 
over} $\{T_j\}$. 
It is easy to see that if $\{T_{i,j}\}$ is controlled over $\{T_j\}$ and $\phi:\{
X_{i,j}\}\to\{Y_j\}$ is a Thom map then that condition is satisfied, i.e., $\{T_{
i,j}\}$ is strongly controlled over $\{T_j\}$. 
Moreover, for a Thom map $\phi:\{X_{i,j}\}\to\{Y_j\}$ and a strongly controlled 
subanalytic $C^m$ tube system $\{T_j\}$ for $\{Y_j\}$ there exists a subanalytic 
$C^m$ tube system $\{T_{i,j}\}$ for $\{X_{i,j}\}$ strongly controlled over $\{T_j
\}$ if the underlying map $\phi:X\to Y$ is of class $C^\omega$ [S$_1$, Lemma 
I.1.3$'$]. 
Thom's conjecture is that a proper Thom map has a triangulation, and [S$_2$] 
showed that if a proper subanalytic $C^\omega$ map $\phi:X\to Y$ admits a Whitney 
subanalytic $C^\omega$ stratification $\phi:\{X_{i,j}\}\to\{Y_j\}$, a strongly 
controlled subanalytic $C^m$ tube system $\{T_j\}$ for $\{Y_j\}$ and a 
subanalytic $C^m$ tube system $\{T_{i,j}\}$ for $\{X_{i,j}\}$ strongly controlled 
over $\{T_j\}$, then $\phi$ has a subanalytic triangulation. 
Here the assumption that $\phi$ is proper is too strong to apply to our case. 
We replace it as follows. 

\proclaim{Theorem (Triangulation of a stratified map)}
Let $\phi:X\to Y$ be a subanalytic $C^\omega$ map, $\phi:\{X_j\}\to\{Y_j\}$ its 
Whitney subanalytic $C^\omega$ stratification such that $X_j=\phi^{-1}(Y_j)$ for 
each $j$, $\{T_{Y,j}=(|T_{Y,j}|,\pi_{Y,j},\rho_{Y,j})\}$ a strongly controlled 
subanalytic $C^m$ tube system for $\{Y_j\}$, and $\{T_{X,j}=(|T_{X,j}|,\pi_{X,j},
\rho_{X,j})\}$ a subanalytic $C^m$ tube system for $\{X_j\}$ strongly controlled 
over $\{T_{Y,j}\}$. 
Assume the map $(\pi_{X,j},\phi)|_{X\cap|T_{X,j}|}:X\cap|T_{X,j}|\to X_j\times|T_
{Y,j}|$ is proper for each $j$. 
Then $\phi$ admits a subanalytic triangulation, and if $X$ is a $C^1$ manifold we 
can choose the triangulation $(\tau_X:P_X\to X,\,\tau_Y:P_Y\to Y)$ of $\phi$ so that 
$P_X$ is PL homeomorphic to the PL manifold of a $C^1$ triangulation of $X$. 
\endproclaim
  We can prove this theorem in the same way as in [S$_2$]. 
Recall the proof in [S$_2$]. 
Consider the case where $X_j=\phi^{-1}(Y_j)$ as in the above theorem, and choose the 
set of indexes so that $\dim Y_j\ge\dim Y_{j+1}$. 
By induction we assume there exists a subanalytic triangulation $(\tau_{X,j}:P_{X,j}
\to\hat X_j\overset\text{def}\to=X-\cup_{i\ge j}\{x\in|T_{X,i}||\,\rho_{Y,i}\circ
\tilde\phi(x)\le\epsilon_i\circ\phi_i\circ\pi_{X,i}(x)\},\,\tau_{Y,j}:P_{Y,j}\to
\hat Y_j\overset\text{def}\to=Y-\cup_{i\ge j}\{y\in|T_{Y,i}||\,\rho_{Y,i}(y)\le
\epsilon_i\circ\pi_{Y,i}(y)\})$ of $\phi|_{\hat X_j}:\hat X_j\to\hat Y_j$ for some 
$j\in\N$, where $\epsilon_i$ are positive subanalytic $C^m$ functions on $Y_i$ so small 
that $\{\epsilon_i\}$ satisfies the condition of a {\it removal data} at page 5 in 
[S$_1$]. 
Then we extend ``canonically'' $\tau_{Y,j}$ to a subanalytic triangulation $\tau_{Y,
j+1}:P_{Y,j+1}\to\hat Y_{j+1}$ by $\rho_{Y,j+1},\ \pi_{Y,j+1}$ and the strong 
controlledness of $\{T_{Y,j+1}\}$. 
Since the extension is canonical we can lift in the same canonical way $\tau_{Y,j+1}$ 
to a subanalytic triangulation $\tau_{X,j+1}$ of $\hat X_{j+1}$ so that $\tau_{X,j+1}$ 
is an extension of $\tau_{X,j}$ and $(\tau_{X,j+1},\,\tau_{Y,j+1})$ is a subanalytic 
triangulation of $\phi|_{\hat X_{j+1}}$ by the strong controlledness of $\{T_{X,i}\}$ 
over $\{T_{Y,i}\}$ and the properness of $\phi$. 
Here the condition of the properness of $\phi$ is used to assure only that the map $(\pi
_{X,j+1},\phi)|_{X\cap|T_{X,j+1}|}:X\cap|T_{X,j+1}|\to X_{j+1}\times|T_{Y,j+1}|$ is 
proper. 
Hence the above theorem holds true. 
 \head
3. Proof of theorem
\endhead
\demo{Proof of triangulation theorem of $\pi:M\to M/G$}
Set $n=\dim M$. \newline
{\it Reduction to the $C^\omega$ case}. 
Assume $M$ is a proper $C^1$ $G$-manifold. 
Then $M$ equivariantly $C^1$ diffeomorphic to some $C^\omega$ $G$-manifold. 
This was shown in Palais [P$_2$] in the case of compact $G$ and $M$, in 
Matumoto-Shiota [M-S$_1$] in the case of compact $G$ and in Illman [I] in the 
general case. However, [I] used without proof the following theorem of Koszul [K]. 
\newline
{\it Fact}. If $X$ is a differentiable $G$-manifold and the isotropy group $G_x$ at 
$x\in X$ is compact, then there exists a near-slice at $x$ $S$ in $X$ (i.e., $x\in S
\subset X,\, G_xS=S$ and there exists a local cross-section $\chi:U\to G$ in $G/G_x$ 
such that the map $U\times S\ni(u,s)\to\chi(u)s\in X$ is a homeomorphism onto an 
open neighborhood of $x$ in $X$). \par
  A problem is that the proof in [K] and [P$_1$] works for $X$ of class $C^k$, $k>1$. 
Hence we prove Fact in the $C^1$ case. 
Regard $X$ as a $C^1$ $G_x$-manifold. 
Note $G_x$ is a compact Lie group, which is a $C^\omega$ submanifold of $G$. 
Then by [M-S$_1$, Theorem 1.3] there exists a $C^\omega$ $G_x$-manifold $X^*$ and a 
$C^1$ $G_x$-equivariant diffeomorphism $f:X\to X^*$. 
Set $f(x)=x^*$ and $f(Gx)=Y$, which is a $C^1$ submanifold of $X^*$ and contains 
$x^*$. 
Choose a $C^\omega$ Riemannian metric for $X^*$ invariant under $G_x$. 
Let $\epsilon>0$ be a small number, define a $C^\omega$ submanifold $S^*$ of $X^*$ 
to be the union of geodesic segments of length$\,<\epsilon$ starting from $x^*$ in 
a direction orthogonal to $T_{x^*}Y$---the tangent space of $Y$ at $x^*$. 
Then $S=f^{-1}(S^*)$ satisfies the requirements. \par
  If $M$ is a subanalytic $C^1$ $G$-manifold, the above proof and the proof in [I] 
work in the subanalytic category and hence $M$ is equivariantly subanalytically 
$C^1$ diffeomorphic to a subanalytic $C^\omega$ $G$-manifold, which is denoted by 
$M'$. (Here $M$ and $M'$ are closedly imbedded in some Euclidean space.) 
We can replace $M$ with $M'$ for the following reason. 
Let $f:M\to M'$ and $\overline f:M/G\to M'/G$ be the $C^1$ diffeomorphism and the 
induced homeomorphism, assume the theorem is 
proved for $M'$ and let $(\tau':P'\to M',\,\sigma':L'\to M'/G)$ be a resulting 
triangulation of the natural map $\pi':M'\to M'/G$. 
Then $(f^{-1}\circ\tau':P'\to M,\,\overline f^{-1}\circ\sigma':L'\to M/G)$ is the 
required triangulation of $\pi:M\to M/G$. 
For that we only need to see $f^{-1}\circ\tau'$ is subanalytic when $M$ is subanalytic. 
It is possible by the fact that the composite of 2 subanalytic maps is subanalytic 
if the first source space is closedly imbedded in a Euclidean space. 
Hence from now we assume $M$ is a proper $C^\omega$ $G$-manifold. 
Let it be closedly imbedded in $\R^{2n+1}$. 
Note $M$ is a subanalytic $G$-manifold even if the original $M$ is not subanalytic. \par
  By [M-S$_2$, Theorem 3.3] there exists a $G$-invariant subanalytic map $p:M\to \R^{
2n+1}$, where $G$ acts on $\R^{2n+1}$ trivially, such that $p(M)$ is closed and 
subanalytic in $\R^{2n+1}$ and the induced map $\overline p:M/G\to p(M)$ is a 
homeomorphism. 
Set $X=\graph p$ and $Y=p(M)$, and let $\phi:X\to Y$ denote the projection, which is 
a subanalytic $C^\omega$ map. 
Note also the action $G\times X\to X$ is of class subanalytic $C^\omega$ because the 
action is $G\times M\times\R^{2n+1}\supset G\times X\ni(g,x,p(x))\to(g x,p(x))\in X$, 
which is the restriction to $G\times X$ of the projection $G\times M\times\R^{2n+1}
\to M\times\R^{2n+1}$. 
Hence the action is extendable to an analytic map $G\times M\times\R^{2n+1}\to M\times
\R^{2n+1}$. 
We will construct a subanalytic triangulation $(\tau:P\to M,\,\sigma:L\to M/G)$ of 
$\phi$. 
Then $L$ is unique by [M-S$_2$, Corollary 3.5]. \par
  There exists the canonical Whitney subanalytic $C^\omega$ stratification of $\phi$ 
for the following reason though $\phi$ does not satisfy the boundedness condition 
$(*)$ at the definition of a stratification of a map. 
Set $X_1=\Reg \phi|_{\Reg X}$. 
Then $X_1$ is a subanalytic $C^\omega$ manifold of dimension $n$, $\dim (X-X_1)<n$, 
and $X_1$ is $G$-invariant because for each $(g,x)\in G\times X$ and for an open smooth 
neighborhood $U$ of $x$ in $X$, the map $U\ni y\to g y\in gU$ is a $C^\omega$ 
diffeomorphism, hence $gU$ is smooth and because $\phi(x)=\phi(g^{-1}x)$ for $x\in gU$. 
Hence $\phi^{-1}(\phi(X_1))=X_1$, $\phi(X_1)$ is a subanalytic $C^\omega$ manifold of 
various local dimension by the same reason and by the fact that $\phi(B)$ is 
subanalytic for each bounded subanalytic set $B\subset X\subset\R^{2n+1}\times\R^{2n+1
}$, and $\phi|_{X_1}$ is a $C^\omega$ submersion onto $\phi(X_1)$. 
Here a $C^\omega$ manifold of various local dimension means a set each of whose 
connected components is a $C^\omega$ manifold. 
For each $j,k=0,...,n$, let $Y_{1,j,k}$ denote the union of connected components $C$ 
of $Y_1\overset\text{def}\to=\phi(X_1)$ of dimension $k$ and such that $\dim\phi^{-1}
(C)=j$. 
Set $X_{1,j,k}=\phi^{-1}(Y_{1,j,k})$. 
Then $\phi|_{X_1}:\{X_{1,j,k}\}\to\{Y_{1,j,k}\}$ is a subanalytic $C^\omega$ 
stratification of $\phi|_{X_1}:X_1\to Y_1$. 
Apply the same arguments to $\phi|_{X-X_1}:X-X_1\to Y-Y_1$. 
Then we have a subanalytic $C^\omega$ manifold $X_2$ in $X-X_1$ and a subanalytic 
$C^\omega$ stratification $\phi|_{X_2}:\{X_{2,j,k}\}\to\{Y_{2,j,k}\}$ of $\phi|_{X_2}
:X_2\to\phi(X_2)$ such that $\phi^{-1}(\phi(X_2))=X_2$ and $\dim(X-X_1-X_2)<\dim(X-X_1)
$. 
We require, moreover, $\{X_{1,j,k},X_{2,j,k}\}$ and $\{Y_{1,j,k},Y_{2,j,k}\}$ to be 
Whitney stratification. 
That is possible because for each $X_{2,j,k}$, the subsets of $X_{2,j,k}$ where some of 
$X_{1,j',k'}$ and $X_{2,j,k}$ do not satisfy the Whitney condition and its closure in 
$X_{2,j,k}$ are $G$-invariant. 
In this way we obtain the canonical Whitney subanalytic $C^\omega$ stratification $\phi
:\{X_{i,j,k}\}\to\{Y_{i,j,k}\}$ of $\phi:X\to Y$, which we write simply $\phi:\{X_j\}
\to\{Y_j\}$. \par
  We see $\phi:\{X_j\}\to\{Y_j\}$ is a Thom map as follows. 
Let $X_j$ and $X_{j'}$ be such that $\overline{X_j}\cap X_{j'}\not=\emptyset$, and let 
$\{a_k\}$ be a sequence in $X_j$ convergent to a point $b$ of $X_{j'}$ such that $\{T_{
a_k}(\phi|_{X_j})^{-1}(\phi(a_k))\}$ converges to a space $T\subset\R^{2n+1}\times\R^{2
n+1}$. 
Write $a_k$ as $(a_k',p(a_k'))\in M\times\R^{2n+1}$. 
Then $(\phi|_{X_j})^{-1}(\phi(a_k))=\phi^{-1}(\phi(a_k))=Ga_k$ since $X_j$ is 
$G$-invariant, and $Ga_k=G a_k'\times\{p(a_k')\}$. 
Hence $T=\lim_{k\to\infty}T_{a'_k}Ga'_k\times\{0\}$. 
Clearly $\lim_{k\to\infty}T_{a'_k}Ga'_k\supset T_{b'}Gb'$ where $b=(b',p(b'))\in M\times
\R^{2n+1}$. 
Therefore, $T_b(\phi|_{X_{j'}})^{-1}(\phi(b))\subset T$. \par
  As noted at the definition of a strongly controlled tube system, there exists a 
strongly controlled subanalytic $C^{2+n}$ tube system $\{T_{Y,j}=(|T_{Y,j}|,\pi_{Y,j},
\rho_{Y,j})\}$ for $\{Y_j\}$. 
Hence by triangulation theorem of a stratified map it remains only to find a subanalytic 
$C^2$ tube system $\{T_{X,j}=(|T_{X,j}|,\pi_{X,j},\rho_{X,j})\}$ for $\{X_j\}$ 
controlled over $\{T_{Y,j}\}$ such that the map $(\pi_{X,j},\phi)|_{X\cap|T_{X,j}|}:X
\cap|T_{X,j}|\to X_j\times|T_{Y,j}|$ is proper for each $j$ since $\phi:\{X_j\}\to\{Y_j
\}$ is a Thom map. 
Moreover, it suffices to define $\pi_{X,j}$ on only $|T_{X,j}|\cap M\times\R^{2n+1}$ 
because if we let $\tilde\phi$ be the projection of $\R^{2n+1}\times\R^{2n+1}$ to the 
latter factor and if $\pi_{X,j}$ on $|T_{X,j}|\cap M\times\R^{2n+1}$ are given so that 
the conditions of controlledness are satisfied then $\pi_{X,j}\circ(q,\id)$ fulfills 
the requirements, where $q$ is the orthogonal projection of a tubular neighborhood of 
$M$ in $\R^{2n+1}$. 
In order to understand the problem of the construction of $\{T_{X,j}\}$ we consider the 
following easy case\,:\par
  {\it Case of compact $G$}. 
Remember that $\phi:X\to Y$ is proper. 
We know there exists a subanalytic $C^2$ tube system $\{T_{X,j}=(|T_{X,j}|,\pi_{X,j},
\rho_{X,j})\}$ for $\{X_j\}$ controlled over $\{T_{Y,j}\}$. 
We need to shrink $|T_{X,j}|$ and $|T_{Y,j}|$ so that the above properness condition is 
satisfied. 
First we can assume $\overline{\phi(X\cap|T_{X,j}|)}\subset|T_{Y,j}|$ and $\pi^{-1}_{Y,
j}(y)$ is bounded for each $y\in Y_j$. 
For each $y\in Y_j$, set 
$$
  \chi(y)=\min\{\rho_{Y,j}\circ\phi(x)\,|\,x\in X\cap\bdry|T_{X,j}|,\pi_{Y,j}\circ\phi(x)
=y\}. $$
Then there exists a positive subanalytic $C^0$ function $\chi'$ on $Y_j$ such that 
$\chi'<\chi$ since $\phi:X\to Y$ is proper. 
Shrink $|T_{Y,j}|$ and $|T_{X,j}|$ to $\{y\in|T_{Y,j}|\,|\,\rho_{Y,j}(y)<\chi'\circ\pi_{Y,
j}(y)\}$ and $|T_{X,j}|\cap\R^{2n+1}\times($new$\,|T_{Y,j}|)$ respectively. 
Then $\phi(X\cap|T_{X,j}|)=Y\cap|T_{Y,j}|$ and $X\cap|T_{X,j}|$ is $G$-invariant. 
Hence $\phi|_{X\cap|T_{X,j}|}:X\cap|T_{X,j}|\to |T_{Y,j}|$ is proper since $\phi:X\to 
Y$ is so. 
It follows that $(\pi_{X,j},\phi)|_{X\cap|T_{X,j}|}:X\cap|T_{X,j}|\to X_j\times|T_{Y,j}
|$ is proper. 
Thus triangulation theorem of $\pi$ follows from triangulation theorem of a stratified 
map. \par
  The above arguments show when $G$ is compact, if we shrink $|T_{Y,j}|$ then we can 
choose enough large domains $|T_{X,j}|$ in comparison. 
However, this is not the case in general. 
In general case we will define $|T_{X,j}|$'s on only slices because isotropy groups are 
compact, and extend it globally. 
For extension we need an additional condition\,:
$$
\pi_{X,j}(g x)=g\pi_{X,j}(x)\quad\text{for}\ (g,x)\in G\times(X\cap|T_{X,j}|). \tag"$(**)
_j$"
$$
This condition can be satisfied in the following case, and we can reduce the case of 
compact $G$ to the following case. \par
  {\it Case where $G$ is a compact subgroup of the orthogonal group $\bold O(2n+1)$, $M
\subset\R^{2n+1}$ and $G$ operates orthogonally on $M$.} 
We construct $\{T_{X,j}\}$ by double induction. 
Choose the set of indexes so that $\dim X_j\le\dim X_{j+1}$. 
Then (sc2) and hence $\rho_{X,j}$ are not necessary since $\phi^{-1}(\phi(X_j))=X_j$, and 
the required conditions are $(**)_j$,
$$
\gather
\phi\circ\pi_{X,j}=\pi_{Y,j}\circ\tilde\phi\quad\text{on}\ |T_{X,j}|,\tag"(sc1)$_j$"\\
\pi_{X,j}\circ\pi_{X,k}=\pi_{X,j}\quad\text{on}\ |T_{X,j}|\cap|T_{X,k}|\ \text{for}\ 
j<k.\tag"(sc3)$_{j,k}$"
\endgather
$$
We require also $\pi_{X,j}$ to be of class subanalytic $C^*\overset\text{def}\to=C^{2+n-
\dim X_j}$ for a technical reason. \par
  Assume for some $k$, $\{T_{X,j}\}_{j<k}$ is given so that (sc1)$_j$, (sc3)$_{l,j}$ and 
$(**)_j$ for any $l<j<k$ are satisfied. 
First, define $\pi_{X,k}(x)$ for $x\in|T_{X,k}|$---a small open neighborhood of $X_k$ 
in $\R^{2n+1}\times\R^{2n+1}$---to be the orthogonal projection of $x$ to the subanalytic 
$C^\omega$ manifold $\phi^{-1}(\pi_{Y,k}\circ\tilde\phi(x))$. 
Then $\pi_{X,k}$ is of class subanalytic $C^{2+n}$, (sc1)$_k$ is satisfied and, moreover, 
so is $(**)_k$ since $G$ acts orthogonally on $M\times\R^{2n+1}$. 
By downward induction we shrink $|T_{X,j}|,\ j\le k$, and modify $\pi_{X,k}$ so that 
(sc3)$_{l,k}$ for $l<k$ are satisfied. 
Assume for some $l<k$ we have shrunk $|T_{X,j}|,\ l<j<k,$ and have modified $\pi_{X,k}$ 
so that (sc3)$_{j,k}$, $l<j<k$, are satisfied and $\pi_{X,k}$ is now of class subanalytic 
$C^*$. 
(In the following arguments we need to shrink $|T_{X,j}|,\,l\le j\le k$, many times. 
However, we do not mention it because it is clear when we need to do.) 
For each $x\in|T_{X,k}|\cap|T_{X,l}|$, the set 
$$
X_{k,l,x}=\pi_{X,l}^{-1}(\pi_{X,l}(x))\cap\phi^{-1}(\pi_{Y,k}\circ\tilde\phi(x))
$$
is a subanalytic $C^{2+n-\dim X_l}$ submanifold of $X_k\cap|T_{X,l}|$ since $\phi$ is a 
Thom map, $\{X_{k,l,x}\,|\,x\in|T_{X,k}|\cap|T_{X,l}|\}$ is a subanalytic $C^{2+n-\dim X
_l}$ foliation of $X_k\cap|T_{X,l}|$, and $x\in X_{k,l,x}$ if $x\in X_k\cap|T_{X,l}|$. 
Let $p_{k,l,x}$ denote the orthogonal projection to $X_{k,l,x}$ of its small 
neighborhood, and set 
$$
\pi_{X,k,l}(x)=p_{k,l,x}(\pi_{X,k}(x))\quad\text{for}\ x\in|T_{X,k}|\cap|T_{X,l}|. 
$$
Then $\pi_{X,k,l}:|T_{X,k}|\cap|T_{X,l}|\to X_k\cap|T_{X,l}|$ is a subanalytic $C^*$ 
submersion, (sc1)$_k$ and (sc3)$_{l,k}$ for $\pi_{X,k,l}$ are clear, and $(**)_k$ for 
$\pi_{X,k,l}$ holds because for $(g,x)\in G\times(X\cap|T_{X,k}|\cap|T_{X,l}|)$ 
$$
\split
\pi_{X,k,l}(g x)=p_{k,l,g x}&(\pi_{X,k}(g x))\\
=p_{k,l,g x}(g\pi_{X,k}(x))&\qquad\text{by}\ (**)_k\ \text{for}\ \pi_{X,k}\\
=g p_{k,l,x}(\pi_{X,k}(x))\ &\qquad\text{since}\ \ p_{k,l,g x}(g y)=g p_{k,l,x}(y)\ \text{for}\ 
y\ \text{near}\ X_{k,l,x}.
\endsplit
$$
Moreover, for any $j$ with $l<j<k$ we have\,:
$$
\pi_{X,k,l}=\pi_{X,k}\quad\text{on}\ |T_{X,k}|\cap|T_{X,l}|\cap|T_{X,j}|\tag$***$
$$
as follows. 
Let $x\in |T_{X,k}|\cap|T_{X,l}|\cap|T_{X,j}|$. 
Then 
$$
\gather
\pi_{X,l}(x)\overset\text{by (sc3)}_{l,j}\to=\pi_{X,l}\circ\pi_{X,j}(x)
\overset\text{by (sc3)}_{j,k}\to=\\
\pi_{X,l}\circ\pi_{X,j}\circ\pi_{X,k}(x)\overset\text{by (sc3)}_{l,j}\to
=\pi_{X,l}\circ\pi_{X,k}(x),\\
\pi_{X,k}(x)\in\pi_{X,l}^{-1}(\pi_{X,l}(x)).\tag"hence" 
\endgather
$$
On the other hand, by (sc1)$_k$
$$
\gather
\pi_{X,k}(x)\in\phi^{-1}(\pi_{Y,k}\circ\tilde\phi(x)). \\
\pi_{X,k}(x)\in X_{k,l,x},\tag"Consequently,"
\endgather 
$$
which proves $(***)$ by definition of $\pi_{X,k,l}$. 
Next we paste $\pi_{X,k}$ and $\pi_{X,k,l}$ by a partition of unity. 
Let $\xi$ be a subanalytic $C^{2+n}$ function on $|T_{Y,l}|$, not $|T_{X,l}|$, such that 
$\xi\ge 0$, $\xi=0$ outside of a small neighborhood of $Y_l$ in $|T_{Y,l}|$, and $\xi=1$ 
on a smaller one. 
Set
$$
\pi'_{X,k}\!(x)\!=\!
\cases
\!\pi_{X,k}(x)\ &\text{for}\,x\!\in|T_{X,k}|-|T_{X,l}|\\
\! p_{k,x}(\xi\circ\tilde\phi(x)\pi_{X,k,l}(x)+(1-\xi\circ\tilde\phi(x))\pi_{X,k}(x))\,&
\text{for}\,x\!\in|T_{X,k}|\cap|T_{X,l}|,
\endcases
$$
where $p_{k,x}$ denotes the orthogonal projection to $\phi^{-1}(\pi_{Y,k}\circ\tilde\phi
(x))$ of its small neighborhood. 
Then $\pi'_{X,k}$ is the required modification of $\pi_{X,k}$. 
Indeed, $\pi'_{X,k}:|T_{X,k}|\to X_k$ is a subanalytic $C^*$ submersion, $\pi'_{X,k}$ 
satisfies clearly (sc1)$_k$ and (sc3)$_{l,k}$, and (sc3)$_{j,k}$ for $l<j<k$ follow from 
$$
\pi'_{X,k}=\pi_{X,k}\quad\text{on}\ |T_{X,k}|\cap|T_{X,j}|,
$$
which is a trivial consequence of $(***)$. 
Lastly, $(**)_k$ for $\pi'_{X,k}$ holds by the same reason as above and by the fact $\xi
\circ\tilde\phi=\const$ on $Gx$ for each $x\in X$. 
Thus there exists $\{T_{X,j}\}$ satisfying $(**)_j$. \par
  {\it Case of non-compact $G$.} 
For each point $a\in M$ there exists a linear slice $S$ at $a$, i.e., a  $C^\omega$ 
submanifold of $M$ such that $S$ contains $a$ and is $G_a$-invariant, the map $G\times_{
G_a}S\ni(g,s)\to g s\in M$ is a $C^\omega$ diffeomorphism onto an open neighborhood of 
$a$ in $M$ and $S$ is $G_a$-equivariantly $C^\omega$ diffeomorphic to a Euclidean space 
where $G_a$ acts orthogonally (see [B]). 
Here $G\times_{G_a}S$ is the quotient space of $G\times S$ under the equivalence relation 
$(g g',s)\sim(g,g's)$ for $(g,g',s)\in G\times G_a\times S$ and we can choose bounded and 
subanalytic $S$. 
Let $a_\alpha,S_\alpha,\alpha\in A$, be a finite or countable number of points of $M$ 
and bounded subanalytic linear slices at $a_\alpha$ such that $\{p(S_\alpha)\}_{
\alpha\in A}$ is a locally finite covering of $Y$. 
The above construction of $\phi:\{X_j\}\to\{Y_j\}$ works with an additional condition 
of compatibility. 
Hence we assume $\phi:\{X_j\}\to\{Y_j\}$ is compatible with $\{X\cap GS_\alpha\times\R^{
2n+1}\}$ and $\{p(S_\alpha)\}$, i.e., each of $X\cap GS_\alpha\times\R^{2n+1}$ or $p(S_
\alpha)$ is a union of some connected components of $X_j$ or $Y_j$ respectively, and 
then for simplicity of notation, each is a union of $X_j$ or $Y_j$ which is possible 
because each connected component of a subanalytic set is subanalytic. 
For each $j$, let $\alpha_j\in A$ be such that $X_j\subset GS_{\alpha_j}\times\R^{2n+1}$. 
Here we can assume, moreover, $\overline{X_j}\subset GS_{\alpha_j}\times\R^{2n+1}$ without 
loss of generality. 
Set $a_j=a_{\alpha_j},\,S_j=S_{\alpha_j},\,G_j=G_{a_j}$ and $Z_j=\graph p|_{S_j}$ 
for each j. 
Choose the set of indexes as in the last case. 
We will construct $\{T_{X,j}\}$ by double induction in the same way. 
Assume for some $k$ we are given $\{T_{X,j}\}_{j<k}$ of class subanalytic $C^*$ such 
that (sc1)$_j$,\,(sc3)$_{l,j}$ and $(**)_j$ for any $l<j<k$ are satisfied and the map 
$(\pi_{X,j},\phi)|_{X\cap|T_{X,j}|}:X\cap|T_{X,j}|\to X_j\times|T_{Y,j}|$ is proper. \par
  We construct $T_{X,k}$ as follows. 
Here also we do not mention shrinking $|T_{X,j}|$ and $|T_{Y,j}|$ each time, though we 
need to keep the properness condition on $(\pi_{X,j},\phi)\allowmathbreak|_{X\cap|T_{X,
j}|}:X\cap|T_{X,j}|\to X_j\times|T_{Y,j}|$. 
First define a subanalytic $C^{1+n}$ tube $T_{Z,k}=(|T_{Z,k}|,\pi_{Z,k},\rho_{Z,k})$ at 
$Z_k\cap X_k$ in $S_k\times\R^{2n+1}$ by $\pi_{Z,k}(z)=$\,the orthogonal projection of 
$z$ to the $C^\omega$ manifold $Z_k\cap\phi^{-1}(\pi_{Y,k}\circ\phi(z))$, where $S_k$ is 
regarded as a Euclidean space and $G_k$ acts orthogonally there and trivially on $\R^{2n
+1}$. 
Then (sc1)$_k$ is satisfied and we have 
$$
\pi_{Z,k}(g z)=g\pi_{Z,k}(z)\quad\text{for}\ (g,z)\in G_k\times(Z_k\cap|T_{Z,k}|)
\tag"$(**)_{Z,k}$"
$$
since $G_k$ acts orthogonally on $Z_k$. 
As in the case of compact $G$, choose $|T_{Y,k}|$ and $|T_{Z,k}|$ so that the map $(\pi_
{Z,k},\phi)|_{Z_k\cap|T_{Z,k}|}:Z_k\cap|T_{Z,k}|\to(Z_k\cap X_k)\times|T_{Y,k}|$ is proper, 
which is possible because $\{g\in G\,|\, gZ_k=Z_k\}=G_k$ and $G_k$ is compact. 
Next we extend $T_{Z,k}$ to a subanalytic $C^{1+n}$ tube $T_{X,k}$. 
For that it suffices to define $\pi_{X,k}$ on $|T_{X,k}|\cap M\times\R^{2n+1}$. 
Choose $|T_{X,k}|$ so that 
$$
|T_{X,k}|\cap M\times\R^{2n+1}=G|T_{Z,k}|,
$$
and set
$$
\pi_{X,k}(g z)=g\pi_{Z,k}(z)\quad\text{for}\ (g,z)\in G\times|T_{Z,k}|. 
$$
First of all $\pi_{X,k}$ is then well-defined. 
Indeed, if $g z=g' z'$ for $(g,z),(g',z')\in G\times|T_{Z,k}|$ then $s$ and $s'$ in $S_k$ 
with $z=(s,p(s))$ and $z'=(s',p(s'))$ satisfy $s=g^{-1}g' s'$ and hence $g^{-1}g'\in G_k$ 
since $G\times_{G_k}S_k\to M$ is an imbedding. 
Therefore, 
$$
g'\pi_{Z,k}(z')=g(g^{-1}g'\pi_{Z,k}(z'))\overset\text{by $(**)_{Z,k}$}\to=g\pi_{Z,k}(g^{-1}g' 
z')=g\pi_{Z,k}(z). 
$$
Clearly $\pi_{X,k}:|T_{X,k}|\cap M\times\R^{2n+1}\to X_k$ is a subanalytic $C^{1+n}$ 
submersion. 
Next
$$
\gather
\phi\circ\pi_{X,k}(g z)=\phi(g\pi_{Z,k}(z))=\phi\circ\pi_{Z,k}(z)\overset\text{by (sc1)$_k$ 
for }\pi_{Z,k}\to=\tag"(sc1)$_k$"\\
\qquad\qquad\pi_{Y,k}\circ\phi(z)=\pi_{Y,k}\circ\phi(g z)\quad\text{for}\ (g,z)\in G\times|T
_{Z,k}|,\\
\pi_{X,k}(g g' z)=g g'\pi_{Z,k}(z)=g\pi_{X,k}(g' z)\qquad\qquad\qquad\tag"$(**)_{X,k}$"\\
\qquad\qquad\qquad\text{for}\ (g,g',z)\in G^2\times|T_{Z,k}|.
\endgather
$$
Lastly, the map $(\pi_{X,k},\phi)|_{X\cap|T_{X,k}|}:X\cap|T_{X,k}|\to X_k\times|T_{Y,k}|$ is 
proper as follows. 
The map $G\times_{G_k}(Z_k\cap|T_{Z,k}|)\ni(g,z)\to g z\in X\cap|T_{X,k}|$ is a diffeomorphism. 
Hence by definition of $\pi_{X,k}$ we can regard $\pi_{X,k}|_{X\cap|T_{X,k}|}$ as the map $G
\times_{G_k}(Z_k\cap|T_{Z,k}|)\ni(g,z)\to (g,\pi_{Z,k}(z))\in G\times_{G_k}(Z_k\cap X_k)$ and 
hence $(\pi_{X,k},\phi)|_{X\cap|T_{X,k}|}$ as the map $G\times_{G_k}(Z_k\cap|T_{Z,k}|)\ni(g,z)
\to (g,\pi_{Z,k}(z),\phi(z))\in G\times_{G_k}((Z_k\cap X_k)\times|T_{Y,k}|)$. 
Therefore, properness follows from that of the map $(\pi_{Z,k},\phi)|_{Z_k\cap|T_{Z,k}|}:Z_k
\cap|T_{Z,k}|\to (Z_k\cap X_k)\times|T_{Y,k}|$ because for any $G_k$-equivariant proper map 
$H_1\to H_2$ between $G_k$-spaces, the map $G\times_{G_k}H_1\to G\times_{G_k}H_2$ is proper. 
\par
  We need to modify $\pi_{X,k}$ so that (sc3)$_{j,k}$ holds for any $j<k$. 
By downward induction we assume for some $l<k$, (sc3)$_{j,k}$ for any $j$ with $l<j<k$ is 
satisfied and $\pi_{X,k}$ is of class subanalytic $C^*$. 
If we modify $\pi_{X,k}$ as in the last case so that (sc3)$_{l,k}$ is satisfied, $|T_{X,k}|$ 
may become too small for the map $(\pi_{X,k},\phi)|_{X\cap|T_{X,k}|}:X\cap|T_{X,k}|\to X_k
\times|T_{Y,k}|$ to be proper. 
However, we can not carry out the modification in $Z_k$---as in the last case---because $S_l
\cap p^{-1}(|T_{Y,l}|\cap|T_{Y,k}|)$ and $S_k\cap p^{-1}(|T_{Y,l}|\cap|T_{Y,k}|)$ may be 
different, though we assume $X_l\subset\overline{X_k}-X_k$ without loss of generality. 
We need a revision to the orthogonal direction to $S_k$ in $M$. 
Regard $M$ as a subanalytic $C^\omega$ $G_k$-manifold, give it a $C^\omega$ $G_k$-invariant 
Riemannian metric [B, Theorem 2.1] and lift trivially the metric to $M\times\R^{2n+1}$. 
Since $Z_k$ is bounded, for each $z\in Z_k\cap|T_{X,k}|\cap|T_{X,l}|$, the set 
$$
\tilde Z_{k,l,z}=\pi_{X,l}^{-1}(\pi_{X,l}(z))\cap\phi^{-1}(\pi_{Y,k}\circ\phi(z))
$$
is a subanalytic $C^{2+n-\dim X_l}$ submanifold of $X_k\cap|T_{X,l}|$ such that $z\in\tilde Z
_{k,l,z}$ if $z\in Z_k\cap X_k\cap|T_{X,l}|$. 
Let $p_{k,l,z}$ denote the projection to $\tilde Z_{k,l,z}$ of its small neighborhood such 
that $p^{-1}_{k,l,z}(z')$ for each $z'\in\tilde Z_{k,l,z}$ is a geodesic curve in the 
neighborhood and orthogonal to $\tilde Z_{k,l,z}$ at $z'$, and set 
$$
\pi_{Z,k,l}(z)=p_{k,l,z}(\pi_{X,k}(z))\quad\text{for}\ z\in Z_k\cap|T_{X,k}|\cap|T_{X,l}|.
$$
Then $\pi_{Z,k,l}:Z_k\cap|T_{X,k}|\cap|T_{X,l}|\to X_k$ is a well-defined subanalytic $C^*$ 
map such that $\pi_{Z,k,l}=\id$ on $Z_k\cap X_k\cap|T_{X,l}|$, (sc1)$_k$ and (sc3)$_{l,k}$ 
with $\pi_{X,l}$ are clear, and $(**)_{Z,k}$ on $G_k\times(Z_k\cap|T_{X,k}|\cap|T_{X,l}|)$ 
holds by the same reason as in the last case. 
The following equality also follows in the same way, for each $j$ with $l<j<k$\,
$$
\pi_{Z,k,l}=\pi_{X,k}\quad\text{on}\ Z_k\cap|T_{X,l}|\cap|T_{X,j}|\cap|T_{X,k}|.
$$
Thus $\pi_{Z,k,l}$ on $Z_k\cap|T_{X,k}|\cap|T_{X,l}|$ has the required properties. 
(Note $\Im\pi_{Z,k,l}$ is not necessarily equal to $Z_k\cap X_k$ and $\pi_{Z,k,l}$ is not a 
submersion to $X_k$.) 
It remains to paste it with $\pi_{Z,k}\overset\text{def}\to=\pi_{X,k}|_{Z_k\cap|T_{X,k}|
\cap|T_{X,l}|}$. 
As before let $\xi$ be a subanalytic $C^{2+n}$ function on $|T_{Y,l}|$ such that $\xi\ge0$, 
$\xi=0$ outside of a small neighborhood of $Y_l$ in $|T_{Y,l}|$ and $\xi=1$ on a smaller one. 
Let $(x,x',t)\in(M\times\R^{2n+1})^2\times[0,\,1]$ be such that $x$ and $x'$ are close each 
other. 
Let $\theta(x,x',t)\in M\times\R^{2n+1}$ denote the point in the shortest geodesic curve 
joining $x$ and $x'$ such that the distance between $x$ and $\theta(x,x',t)$ equals the 
product of $t$ and the distance between $x$ and $x'$. 
Set 
$$
\pi'_{Z,k}(z)=
\cases\pi_{Z,k}(z)\quad&\text{for}\ z\in Z_k\cap|T_{X,k}|-|T_{X,l}|\\
p_{Z,k,z}\circ\theta(\pi_{Z,k}(z),\pi_{Z,k,l}(z),\xi\circ\phi(z))\ &\text{for}\ z\in Z_k
\cap|T_{X,k}|\cap|T_{X,l}|,
\endcases
$$
where $p_{Z,k,z}$ denotes the projection to $\phi^{-1}(\pi_{Y,k}\circ\phi(z))$ of its small 
neighborhood with the same properties as $p_{k,l,z}$. 
Then $\pi'_{Z,k}:Z_k\cap|T_{X,k}|\to X_k$ is a subanalytic $C^*$ map, and $\pi'_{Z,k}=\id$ 
on $Z_k\cap X_k$, (sc1)$_k$ for $\pi'_{Z,k}$ on $Z_k\cap|T_{X,k}|$, (sc3)$_{j,k}$ for $\pi'
_{Z,k}$ and $\pi_{X,j},\ l\le j<k$, on $Z_k\cap|T_{X,j}|\cap|T_{X,k}|$ and $(**)_{Z,k}$ for 
$\pi'_{Z,k}$ all hold. 
Extend $\pi'_{Z,k}$ to a subanalytic $C^*$ map $\pi'_{X,k}:|T_{X,k}|\cap M\times\R^{2n+1}
\to X_k$ in the same way as before by 
$$
\pi'_{X,k}(g z)=g\pi'_{Z,k}(z)\quad\text{for}\ (g,z)\in G\times(Z_k\cap|T_{X,k}|).
$$
Then by the same reason, $\pi'_{X,k}$ is a well-defined submersion; (sc1)$_k$ and $(**)_{X,
k}$ are satisfied; the map $(\pi'_{X,k},\phi)|_{X\cap|T_{X,k}|}:X\cap|T_{X,k}|\to X_k
\times|T_{Y,k}|$ is proper; (sc3)$_{j,k}$, $l<j<k$, hold because
$$
\pi'_{X,k}=\pi_{X,k}\quad\text{on}\ |T_{X,k}|\cap|T_{X,j}|\cap M\times\R^{2n+1};
$$
finally
$$
\gather
\pi_{X,l}\circ\pi'_{X,k}(g z)=\pi_{X,l}(g\pi'_{Z,k}(z))\overset\text{by $(**)_{X,l}
$}\to=g\pi_{X,l}\circ\pi'_{Z,k}(z)\tag"(sc3)$_{l,k}$"\\
\overset\text{by (sc3)$_{l,k}$ for $\pi'_{Z,k}$}\to
=g\pi_{X,l}(z)\overset\text{by $(**)_{X,l}$}\to=\pi_{X,l}(g z)\quad\text{for}\ (g,z)\in 
G\times(Z_k\cap|T_{X,k}|\cap|T_{X,l}|).
\endgather
$$
Thus $\pi'_{X,k}$ is the required modification of $\pi_{X,k}$ in the induction process, 
and hence we complete the proof.
\qed
\enddemo

\enddocument